\documentclass[12pt]{amsart}
\usepackage{a4wide}
\usepackage{amssymb, amsmath}
\usepackage{color}
\usepackage{amsthm}
\usepackage{comment}

\newtheorem{thm}{Theorem}[section]
\newtheorem{lem}[thm]{Lemma}

\newtheorem{question}[thm]{Question}

\newcommand{\blokje}{\hfill $\Box$\\}

\begin{document}
\title{Biased partitions of $\mathbb{Z}^n$}
\author{Peter van Hintum}
\thanks{Department of Pure Mathematics and Mathematical Statistics, Centre for Mathematical Sciences, University of Cambridge, Wilberforce road, Cambridge CB3 0WA, United Kingdom. E-mail: pllvanhintum@maths.cam.ac.uk\\
The author was supported by the Cambridge Trust and the UK Engineering and Physical Sciences Research
Council.}
\maketitle

\begin{abstract}
Given a function $f$ on the vertex set of some graph $G$, a scenery, let a simple random walk run over the graph and produce a sequence of values. Is it possible to, with high probability, reconstruct the scenery $f$ from this random sequence?

To show this is impossible for some graphs, Gross and Grupel, in \cite{gross2017}, call a function $f:V\to\{0,1\}$ on the vertex set of a graph $G=(V,E)$ \emph{$p$-biased} if for each vertex $v$ the fraction of neighbours on which $f$ is 1 is exactly $p$. Clearly, two $p$-biased functions are indistinguishable based on their sceneries. Gross and Grupel construct $p$-biased functions on the hypercube $\{0,1\}^n$ and ask for what $p\in[0,1]$ there exist $p$-biased functions on $\mathbb{Z}^n$ and additionally how many there are. We fully answer this question by giving a complete characterization of these values of $p$. We show that $p$-biased functions exist for all $p=c/2n$ with $c\in\{0,\dots,2n\}$ and, in fact, there are uncountably many of them for every $c\in\{1,\dots,2n-1\}$. To this end, we construct uncountably many partitions of $\mathbb{Z}^n$ into $2n$ parts such that every element of $\mathbb{Z}^n$ has exactly one neighbour in each part. This additionally shows that not all sceneries on $\mathbb{Z}^n$ can be reconstructed from a sequence of values on attained on a simple random walk.
\end{abstract}

\section{Background}
For a graph $G=(V,E)$, we call a function $f$ on $V$ a \emph{scenery}. Let $\tilde{X}=(X_n)_{n=0}^\infty$ be a simple random walk on $G$. We associate with $\tilde{X}$ the sequence $(f(X_n))_{n=0}^\infty$ of values attained by the scenery. Is it possible to, with high probability, reconstruct the scenery $f$ from this random sequence? This question has an extensive history, in particular, the case where $G=\mathbb{Z}$ has been studied in many papers, see e.g. \cite{benjamini1996,finucane2014,Howard1996,lindenstrauss1999,matzinger2006} . In \cite{Howard1996}, Howard showed that periodic sceneries on $\mathbb{Z}$ or equivalently sceneries on finite cycles can always be reconstructed. Matzinger and Lember \cite{matzinger2006} extended this result on periodic sceneries to random walks on $\mathbb{Z}$ which include steps of different sizes. In turn, Finuncane, Tamuz and Yaar \cite{finucane2014} refined this idea and extended it to more general Cayley graphs of finite abelian groups.

Focusing on the question for $G=\mathbb{Z}$ (and $G=\mathbb{Z}^2$), Benjamini and Kesten\cite{benjamini1996} showed that almost all sceneries can in fact be distinguished, in the sense that for any given scenery on $G$, it can be distinguished with probability 1 from a scenery chosen randomly in the product measure on all sceneries. Benjamini and, independently, Keane and Den Hollander  conjectured that in fact all pairs of sceneries on $\mathbb{Z}$ are distinguishable (unpublished) \cite{howard1996orthogonality}. This, however, was soon disproved by Lindenstrauss \cite{lindenstrauss1999} who constructed a collection of uncountably many distinct yet indistinguishable sceneries.

In a recent paper, Gross and Grupel \cite{gross2017} showed that for 0-1 functions on the hypercube, i.e. functions of the form $f:\{0,1\}^n\to \{0,1\}$, this reconstruction is not possible in general.  To this end, they defined a $p$-biased function to be a function which takes the value 1 on exactly a fraction $p$ of the neighbours of each vertex in the graph. Note that if there exist two non-isomorphic  $p$-biased functions on a graph, they will be indistinguishable as the sequence of values for both will look like a sequence of independent Bernoulli random variables with succes probability $p$.

Gross and Grupel \cite{gross2017} extended their construction from the hypercube to $\mathbb{Z}^n$ to find  $p$-biased functions for $p=c/2^k$ where $2^k|n$. They  asked for what $p\in [0,1]$ such  $p$-biased functions exist and how many there are for each $p$. Our main aim in this paper is to give a complete characterization of all the values $p\in[0,1]$ for which a  $p$-biased function on $\mathbb{Z}^n$ exists and the exact number for each of those $p$.

As every element of $Z^n$ has $2n$ neighbours, if there is a  $p$-biased function, then $p$ is of the form $p=\frac{c}{2n}$ for some integer $c$. In fact, we will find that for all $c\in\{0,\dots,2n\}$, there exists a  $p$-biased function with $p=\frac{c}{2n}$. Note that the sum of a  $p$-biased function and a  $q$-biased function on disjoint supports is a $p+q$-biased function. Hence, to show that these  $\frac{c}{2n}$-biased functions exist, it suffices to construct $2n$   $\frac{1}{2n}$-biased functions the supports of which partition $\mathbb{Z}^n$.\\
We proceed by showing that for all $n>1$ there are uncountably many non-isomorphic such partitions and uncountably many indistinguishable $p$-biased functions, answering the second question from \cite{gross2017}. Our result extends the result of Lindenstrauss\cite{lindenstrauss1999} that there exist uncountably many indistinguishable sceneries on $\mathbb{Z}$.

\section{Definitions}
As usual, we write $[n]=\{1,\dots, n\}$ and $e_i$ for the $i$th basis vector in the canonical basis of $\mathbb{Z}^n$. We consider $\mathbb{Z}^n$ to be the graph on that  set with edge set $E=\{\{x,x+e_i\}:x\in \mathbb{Z}^n, i\in[n]\}$. Given vertex $x$, $\Gamma(x)$ denotes the neighbourhood, in particular for $x\in \mathbb{Z}^n$; $\Gamma(x)=\{y:y=x\pm e_i \text{ for some } i\in [n]\}$. We identify a 0-1 function with its support, so that we can talk of sets rather than functions.\\
A set $X\subset \mathbb{Z}^n$ is \emph{ $p$-biased} if for all $x\in\mathbb{Z}^n$, we have $|\Gamma(x)\cap X|=2pn$. A partition $\{X_i\}_{i\in[2n]}$ of $\mathbb{Z}^n$ is \emph{biased} if each of its elements is  $\frac{1}{2n}$-biased. 
Hence, as noted in the introduction, if there is a biased partition of $\mathbb{Z}^n$, then there are  $p$-biased function on $\mathbb{Z}^n$ for all $p=c/2n$.\\
To stress that a certain union consists of disjoint sets, we use $\sqcup$ rather than the more common $\cup$.

\section{There is a biased partition of $\mathbb{Z}^n$ for every $n\in\mathbb{N}$}
\begin{thm}\label{main}
For every $n\in\mathbb{N}$, there is a biased partition of  $\mathbb{Z}^n$.
\end{thm}

The proof of this result will consist of a recursive construction. Accordingly, we start by observing that $\mathbb{Z}=\{m: m\equiv 0,1 \mod 4\}\sqcup \{m:m\equiv 2,3\mod 4\}$ is a biased partition. For the recursive construction of the biased partitions, we define a closely related notion. Let $m,n\in\mathbb{N}$ such that $m$ is a multiple of $n$ and let $\{X^i_j\}_{i\in[\frac{m+n}{n}],j\in[2n]}$ be a family of subsets of $\mathbb{Z}^m$. We write $X^i=\bigsqcup_j X^i_j$. We say this family of subsets of $\mathbb{Z}^m$ is \emph{$(m,n)$-filling} if the family partitions $\mathbb{Z}^m$ and if, for each $i\in[\frac{m+n}{n}]$ and $ j\in[2n]$, we have that if $x\in \mathbb{Z}^m\setminus X^i$, then $|\Gamma(x)\cap X^i_j|=1$  and  if $x\in X^i$, then $\Gamma(x)\cap X^i=\emptyset$.\\
 The following lemma is the basis for our recursive construction.

\begin{lem}\label{fillingbp}
If there is a biased partition of $\mathbb{Z}^n$ and there exists an  $(m,n)$-filling family of subsets of $\mathbb{Z}^m$, then there is a biased partition of $\mathbb{Z}^{m+n}$.
\end{lem}
Proof. Let $\{X^i_j\}_{i\in[\frac{m+n}{n}],j\in[2n]}$ be an $(m,n)$-filling family of subsets of $\mathbb{Z}^m$ and let $\{Y_i\}_{i\in[2n]}$ be a biased partition of $\mathbb{Z}^n$ that witnesses $n$ allowing a biased partition. Let, for $i\in [\frac{m+n}{n}]$ and $l\in[2n]$;
\begin{align}\label{Zdef}
Z^i_l=\bigsqcup_j X^i_{j+l}\times Y_j
\end{align}
we claim this is set is  $\frac{1}{2(m+n)}$-biased in $\mathbb{Z}^{m+n}$. For notational convenience, we take $i=1$ and $l=0$.

We claim that every $z=(x,y)\in \mathbb{Z}^{m}\times\mathbb{Z}^n$ has a unique neighbour in $Z=Z^1_0$. If $x\in X^1_j$ for some $j$, then $z$ cannot have a neighbour in $Z$ in the first $m$ coordinates, as $\Gamma(x)\cap X^1=\emptyset$. Fortunately, by definition of $Y_j$, $y$ has a unique neighbour in $Y_j$, say $y'$, which gives $z$ a unique neighbour in the last $n$ coordinates, i.e. $(x,y')$.\\
If, on the other hand, $x\not\in X^1=\bigsqcup_jX^1_j$, then $z$ cannot have a neighbour in $Z$ in the last $n$ coordinates. Since the sets $Y_j$ partition $\mathbb{Z}^n$, we find a unique $j$ such that $y\in Y_j$. By construction, we know that exactly one neighbour of $x$ is in $X^1_j$, say $x'$, so that neighbour gives a unique neighbour of $z$ in $Z$, i.e. $(x',y)$.

Analogously, $Z^i_l$ is a  $\frac{1}{2(m+n)}$-biased set for each $ i\in[\frac{m+n}{n}], l\in [2n]$. Moreover, these sets partition  $\mathbb{Z}^{m+n}$, showing that there is a biased partition of $\mathbb{Z}^{m+n}$. \blokje

For the proof of Theorem 3.1, it remains to find suitable $(m,n)$-filling families.
\begin{lem}\label{timestwo}
Given $n\in\mathbb{N}$, define for $l\in[2],j\in[2n]$ the following subsets of $\mathbb{Z}^n$:
$$X^l_j=
\begin{cases}
\{x\in\mathbb{Z}^n: \sum x_i\equiv l \mod 4,\ \ \ \ \ \ \  \sum ix_i\equiv j \mod n\}\ \ \ \  \text{if $j\leq n$}\\
\{x\in\mathbb{Z}^n: \sum x_i\equiv l +2\mod 4,\ \  \sum ix_i\equiv j \mod n\}\ \ \  \ \text{if $j>n$}
\end{cases}$$
 Then the family $\{X^l_j\}_{l\in[2],j\in[2n]}$ is an $(n,n)$-filling family.
\end{lem}
Proof.  Note that $X^l=\bigsqcup_j X^l_j=\{x:\sum x_i\equiv l \mod 2\}$, so the sets $X_j^l$ partition $\mathbb{Z}^n$. \\
Let $x\in \mathbb{Z}^n$ and $l\in[2]$ be such that $x\in X^l$. Then $x\pm e_i$, the neighbours of $x$, are in $X^{3-l}$ for all $i\in [n]$. It remains to show that these neighbours are all in distinct $X^{3-l}_k$. We find that the neighbour $y=x+e_j$ is such that $\sum y_i\equiv 1+\sum x_i \mod 4$ and $\sum iy_i\equiv j+\sum ix_i\mod n$. For distinct $j\in[n]$ the second sums are clearly distinct modulo $n$. If we compare $y$ to $z=x-e_k$ for some $k\in[n]$, we find that $\sum z_i\equiv-1+\sum x_i\mod 4$, so $y$ and $z$ belong to different parts of the partition. 
\blokje

These lemmas imply that if there is a biased partition of $\mathbb{Z}^n$, then there also is a biased partition of $\mathbb{Z}^{2n}$. Hence, there is a biased partition of $\mathbb{Z}^{2^k}$ for all $k\in\mathbb{N}$. To extend this to the natural numbers with odd prime divisors, we use the following construction.

\begin{lem} \label{epic}
There exists a $(2mn,n)$-filling family of subsets of $\mathbb{Z}^{2mn}$ for all $m,n\in \mathbb{N}$.
\end{lem}
Proof. Define for $l\in[2m+1]$ and $k\in[2n]$;
\begin{align*}
X^l&=\left\{x\in \mathbb{Z}^{2mn}: \sum_{j=1}^{m}\sum_{i=2(j-1)n+1}^{2nj} jx_i \equiv l \mod 2m+1\right\},\\
X^l_k&=\left\{x\in \mathbb{Z}^{2mn}: x\in X^l,\ \ \sum_{i=1}^{2mn} ix_i \equiv k\mod 2n\right\}
\end{align*}
We will show that these sets $X^l_k$  form a $(2mn,n)$-filling family. To this end, note that the sets $X^l$ partition $\mathbb{Z}^{2mn}$ and that   $\{X^l_k\}_{k\in[2n]}$ is a partition of $X^l$ into $2n$ parts. It remains to check that each element of $X^l$ has exactly one neighbour in each of the parts $X^{l'}_k$ for $l'\neq l$ and no neighbour in $X^l$.

Fix some $l\in[2m+1]$. Consider $y\in \mathbb{Z}^{2mn}$. If  $\sum_{j=0}^{m-1}\sum_{i=2jn+1}^{2n(j+1)} jy_i \equiv l \mod 2m+1$, we have $y\in X^l$, and any neighbour of $y$ is not in $X^l$ as changing any coordinate would change this sum.

If $\sum_{j=0}^{m-1}\sum_{i=2jn+1}^{2n(j+1)} jx_i \equiv l-j\mod 2m+1$, for some $j\in [2m]$, there are two options.\\
If $j\in [m]$, then we know that $y+e_i\in X^l$ for every $i\in \{2jn+1,\dots,2(j+1)n\}$. Note that each of these $2n$ vectors is in a different set $X^l_k$. \\
Finally, if $j\in\{m+1,\dots,2m\}$, then let $h=2m+1-j$. Then $y-e_i\in X^l$ for every  $i\in \{2hn+1,\dots,2(h+1)n\}$. Again, each of these $2n$ vectors is in a different set $X^l_j$.
Note that in both cases, these are the only neighbours of $y$ in $X^l$.
\blokje

This is the last ingredient needed for the proof of Theorem \ref{main}.

Proof. As noted, there is a biased partition of $\mathbb{Z}$. In combination with Lemma \ref{timestwo} and Lemma \ref{fillingbp}, this implies that there is a biased partition of $\mathbb{Z}^{2^k}$ for every $k\in \mathbb{Z}_{\geq 0}$. To extend this to all of $\mathbb{N}$, let $n=2^k(2m+1)$ for some $k,m\in\mathbb{Z}_{\geq 0}$. By Lemma \ref{epic}, we can find a $(2m2^k,2^k)$-filling family, which by Lemma \ref{fillingbp} implies that there is a biased partition of $\mathbb{Z}^{2^k+2m2^k}=\mathbb{Z}^n$.\blokje

\section{Counting biased partitions}

In \cite{gross2017}, Gross and Grupel ask, besides a characterisation of $p$ values for which $p$-biased functions exist, for a count of the number of non-isomorphic such  $p$-biased functions. In fact, they provide some finite lower bounds on this number for $p=\frac1n$ and $p=\frac12$, based on an extension of their construction on the hypercube. We find the following complete characterization.
\begin{thm}\label{countingfunctions}
For all $p=c/2n$ with $n>1$ and $c\in\{1,\dots,2n-1\}$, there are $2^{\aleph_0}$ non-isomorphic  $p$-biased functions in $\mathbb{Z}^n$.
\end{thm}

 We say two biased partitions $\{X_i\}_{i\in[2n]}$ and $\{Y_i\}_{i\in[2n]}$ of $\mathbb{Z}^n$ are isomorphic if there exist a graph isomorphism $\phi:\mathbb{Z}^n\to\mathbb{Z}^n$ and permutation $\sigma:[2n]\to[2n]$ such that $\phi(X_i)=Y_{\sigma(i)}$ for all $i\in[2n]$.

\begin{thm}\label{counting}
For all $n>1$ , there are $2^{\aleph_0}$ non-isomorphic biased partitions of $\mathbb{Z}^n$.
\end{thm}

As there are only $2^{\aleph_0}$ subsets of $\mathbb{Z}^n$, the upper bound for so for both Theorems \ref{countingfunctions} and \ref{counting}  are immediate.

In fact, our construction proving Theorems \ref{countingfunctions} and \ref{counting} will be almost identical to the one in the previous section. By introducing a degree of freedom, we produce an uncountable collection of distinct biased partitions. An automorphism $\phi:\mathbb{Z}^n\to\mathbb{Z}^n$ consists of three components; a permutation of the coordinates ($n!$), a reflection of some of the coordinates ($2^n$) and a translation ($|\mathbb{Z}^n|=\aleph_0$). There are only countably many combinations  of these operations and hence there are only countably many automorphisms $\phi:\mathbb{Z}^n\to\mathbb{Z}^n$.

To use the previous construction, we need to have a sense of isomorphism for $(m,n)$-filling families.
We say two $(m,n)$-filling families $\{A^i_j\}_{i\in[\frac{m+n}{n}],j\in[2n]}$ and $\{B^i_j\}_{i\in[\frac{m+n}{n}],j\in[2n]}$ are isomorphic, if there exist a graph isomorphism $\phi:\mathbb{Z}^{m}\to\mathbb{Z}^m$, a permutation $\sigma:[\frac{m+n}{n}]\to[\frac{m+n}{n}]$ and a family of permutations $\tau_i:[2n]\to[2n]$ for all $i\in[\frac{m+n}{n}]$, such that $\phi(A^i_j)=B^{\sigma(i)}_{\tau_{i}(j)}$. Note that by the above observation, each $(m,n)$-filling family is isomorphic to at most countably many other $(m,n)$-filling families.

It is interesting to note the following. Given an $(m,n)$-filling family, we can define a function, taking any element of $\mathbb{Z}^m$ to the the part of the partition that element is in. Two functions arising in such a way from two non-isomorphic $(m,n)$-filling families will be indistinguishable by a simple random walk.

\begin{lem}\label{noniso}
If there is a collection of pairwise non-isomorphic $(m,n)$-filling families of size $2^{\aleph_0}$ and a biased partition of $\mathbb{Z}^n$, then there are $2^{\aleph_0}$ pairwise non-isomorphic biased partitions of $\mathbb{Z}^{m+n}$.
\end{lem}
Proof. In fact, the construction in Lemma \ref{fillingbp} produces such a collection. Let $\{X^i_{j,x}\}_{i\in[\frac{m+n}{n}],j\in[2n]}$, be an uncountable family of pairwise non-isomorphic $(m,n)$-filling families indexed by some $x\in\mathbb{R}$ and let $\{Y_i\}_{i\in [2n]}$ be a biased partition of $\mathbb{Z}^n$.  Let $\{Z^i_{l,x}\}_{i\in[\frac{m+n}{n}],l\in [2n]}$ be the biased partitions as defined in equation (\ref{Zdef}).

As only countably many of the $\{Z^i_{l,x}\}_{i\in[\frac{m+n}{n}],l\in [2n]}$ are isomorphic, it suffices to prove there are $2^{\aleph_0}$ distinct biased partitions.

We claim that all $(m,n)$-filling families $\{X^i_{j,x}\}_{i\in[\frac{m+n}{n}],j\in[2n]}$  produce distinct biased partitions. For a contradiction, suppose $\{Z^i_{l,x}\}_{i\in[\frac{m+n}{n}],l\in [2n]}=\{Z^i_{l,y}\}_{i\in[\frac{m+n}{n}],l\in [2n]}$ for some $x\neq y$. Then we find that for any $i\in[\frac{m+n}{n}]$ and $l\in[2n]$, we can find some $i'\in[\frac{m+n}{n}]$ and $l'\in[2n]$ such that;
$\bigsqcup_j X^{i}_{j+l,x}\times Y_j=\bigsqcup_j X^{i'}_{j+l',y}\times Y_j$. As we know the $Y_j$, we can then find for all $j\in[2n]$; $X^i_{j+l,x}= X^{i'}_{j+l',y}$ and $\{X^i_{j,x}\}_{i\in[\frac{m+n}{n}],j\in[2n]}=\{X^i_{j,y}\}_{i\in[\frac{m+n}{n}],j\in[2n]}$. However, we assumed that the $(m,n)$-filling families were distinct; this contradiction proves the lemma.
\blokje

The construction of Lemma \ref{epic} works in such a way that for any element not in $X^l$ all $2n$ neighbours in $X^l$ lie in the same hyperplane in $\mathbb{Z}^{m+n}$ defined by  $\sum_{j=0}^{m-1}\sum_{i=2jn+1}^{2n(j+1)} jx_i =l+(2m+1)h$ for some $h$. This has the advantage that the construction used to make sure that each set $X^l_k$  contains exactly one of those neighbours need not be the same on distinct planes, giving the following construction.

\begin{lem}\label{2mninfty}
For any  $m,n\in \mathbb{N}$, there are $2^{\aleph_0}$ $(2mn,n)$-filling families.
\end{lem}
Proof. Consider the following slight alteration of the construction in Lemma \ref{epic}.
\begin{align}\label{2mneq}
X^l_{k,f}&=\left\{x\in \mathbb{Z}^{2mn}:\exists h\in \mathbb{Z};  \sum_{j=0}^{m-1}\sum_{i=2jn+1}^{2n(j+1)} jx_i = l+h(2m+1) \text{ and } \sum_{i=1}^{2mn} ix_i \equiv k+f(h)\mod 2n\right\}
\end{align}
where $f:\mathbb{Z}\to[2n]$ is any function. This family is $(2mn,n)$-filling by the proof of Lemma \ref{epic}. There are $\big|[2n]^{\mathbb{Z}}\big|=2^{\aleph_0}$ such functions, and thus such biased partitions. As each of these biased partitions is isomorphic to at most countably many others, we must have $2^{\aleph_0}$ non-isomorphic biased partitions among these.\blokje
 
 Similarly we can extend the construction from Lemma \ref{timestwo}.
 
 \begin{lem}\label{nninfty}
 For $n>1$, there are $2^{\aleph_0}$ $(n,n)$-filling families
 \end{lem}
 Proof. Consider the following construction for $l\in[2]$, $p\in\{0,1\}$, $q\in[n]$ and $f:\mathbb{Z}\to[n]$:
 \begin{align}\label{nneq}
 X^l_{p,q,f}=\left\{x\in \mathbb{Z}^{n}:\exists h\in \mathbb{Z};  \sum_{i=1}^{n} x_i = l+2p+4h \text{ and } \sum_{i=1}^{n} ix_i \equiv q+f(h)\mod n\right\}\end{align}
Writing it in proper form, let $X^l_{k,f}=X^l_{p,q,f}$ with $p=\begin{cases} 0 \text{ if } k\leq n\\ 1 \text{ if } k> n\end{cases}$, and $q\in[n]$ with $q\equiv k\mod n$.
To check that this is in fact a biased partition, note that $X^l_f=\bigsqcup_k X^l_{k,f}=\{x:\sum x_i\equiv l\mod 2\}$, so the $X^l_{k,f}$ form a partition of $\mathbb{Z}^n$ and $\Gamma(X^l_f)\cap X^l_f=\emptyset$. We proceed to check that each element of $\mathbb{Z}^n\setminus X^l_f$ has a unique neighbour in $X^{l}_f$.

Let $x\in \mathbb{Z}^n\setminus X^l_f$, i.e. $\sum x_i\equiv l+1\mod 2$. Write $X^l_{p,f}=\bigsqcup_q X^l_{p,q,f}$, then $x-e_j\in X^l_{p,f}$  and $x+e_j\in X^l_{1-p,f}$ for all $j\in[n]$ for some $p\in\{0,1\}$. Let $y=x+e_j$, then we have $\sum_i y_i=1+\sum_{i} x_i = l+2(1-p)+4h$ for some $h$ not dependent on $j$. Thus, for different $j\in[n]$, we find  $\sum_i i y_i=j+\sum_{i=1}^{n} ix_i \equiv q+f(h)\mod n$ with distinct $q$. Hence, $y=x+e_j\in X^l_{1-p,q,f}$ with distinct  $q$ for  distinct $j$.

Analogously $x-e_j\in X^l_{p,q,f}$ with distinct $q$ for distinct $j$.

As in the proof of Lemma \ref{2mninfty}, we find that for $n>1$, there are $2^{\aleph_0}$ functions $f:\mathbb{Z}\to [n]$ and thus sets $X^l_{p,q,f}$. At most $\aleph_0$ of those can be pairwise isomorphic, so there must be $2^{\aleph_0}$ pairwise non-isomorphic $(n,n)$-filling families.\blokje

What remains is to count the number of non-isomorphic biased partitions of $\mathbb{Z}^2$.

\begin{lem}\label{counting2}
There are $2^{\aleph_0}$ non-isomorphic biased partitions of $\mathbb{Z}^2$.
\end{lem}
Proof. Consider the following set $S=\{x\in\mathbb{Z}^2: x_1+x_2=0\text{ and }x_1\equiv 0 \mod 2\}$. Note that every element of the set $\{x\in\mathbb{Z}^2: x_1+x_2\in\{-1,0,1,2\}\}$, i.e.  four downward diagonals around  the origin, has exactly one neighbour in the set $S+\{0,e_i\}$ for both $i=1$ and $i=2$. 

Given some function $f:\mathbb{Z}\to[2]$, let $X_f=S+\{(2n,2n),(2n,2n)+e_{f(n)}: n\in \mathbb{Z}\}$. Now consider the following biased partition:
$$X_f^k=\begin{cases}
 X_f \text{ if } k=1\\
X_f+(1,-1) \text{ if } k=2\\
X_f+(1,1) \text{ if } k=3\\
x_f+(2,0) \text{ if } k=4
\end{cases}
$$
These 4 sets partition $\mathbb{Z}^2$. Each of the sets is  $\frac{1}{4}$-biased, by the note above. Finally this produces $2^{\aleph_0}$ non-isomorphic biased partitions by the same argument as for the previous two lemmas. \blokje

Proof (of Theorem \ref{counting}).\\
By Lemma \ref{counting2}, we find uncountably many non-isomorphic biased partitions of $\mathbb{Z}^2$, and by Lemmas \ref{2mninfty} and \ref{nninfty} combined with Lemma \ref{noniso}, we find uncountably many non-isomorphic biased partitions of $\mathbb{Z}^n$ for $n>2$. \blokje

Extending the theorem on biased partitions to  $p$-biased functions is not immediate as different biased partitions might give rise to the same  $p$-biased functions. Consider for instance the  $\frac12$-biased functions on $\mathbb{Z}^2$ which has support on $X^1_f\cup X^2_f$ from the proof of Lemma \ref{counting2}; this function is the same for all choices of $f$. We will see that for dimensions bigger than two this is not a problem. We consider the case for $\mathbb{Z}^2$ seperately.
\begin{lem}\label{2function}
For all $p\in\{\frac14,\frac12,\frac34\}$ , there exist infinitely many non-isomorphic  $p$-biased functions in $\mathbb{Z}^2$.
\end{lem}
Proof. For $p=\frac14$ and $p=\frac34$, this follows immediately from Lemma \ref{counting2}, so consider $p=\frac12$. Let $f:\mathbb{Z}\to [2]$ be any function, consider set:
$$X_f=\{x\in\mathbb{Z}^2: x_1\equiv f(x_1+x_2) \mod 2\}$$
Note that for any $x\in\mathbb{Z}^2 $ either $x+e_1$ or $x+e_2$ is in $X_f$, and similarly either $x-e_1$ or $x-e_2$ is in $X_f$. Hence, $X_f$ is  $\frac12$-biased. As there are $2^{\aleph_0}$ different functions $f:\mathbb{Z}\to[2]$, there are as many different  $\frac12$-biased functions, at most countably many of which are isomorphic. \blokje

All that remains is to prove Theorem \ref{countingfunctions}

Proof. Note that for $p=\frac{1}{2n}$ and $p=\frac{2n-1}{2n}$ this follows immediately from Theorem \ref{counting}. The case $n=2$ follows from Lemma \ref{2function}.  Fix some $n>2$ and $p=c/2n$. Use the constructions of equations \ref{2mneq} and \ref{nneq} to produce biased partitions indexed by some function $f$ to feed into the construction in equation \ref{Zdef}. This produces the sets $Z^l_{k,f}$ with $l\in[r]$ and $k\in[q]$ with $r,q$ integers depending on the last step of the recursive construction of biased partitions. That is; $r=2$ and $q=n$ if $n$ is some power of 2 and $r=m+1$ and $q=2^{k+1}$ if $n=2^k(2m+1)$ for some $k\in\mathbb{Z}_{\geq0}$ and $m\in \mathbb{N}$.

Let $I\subset [r]\times[q]$ be such that $|I|=c$ and there is an $l_0\in[r]$ such that $|\{k\in [q]: (l_0,k)\in I\}|$ is either 1 or $q-1$. Let
$$S_f=\bigsqcup_{(l,k)\in I} Z^l_{k,f}$$
We claim that $\{S_f: f\in[r]^\mathbb{Z}\}$ is a family of uncountably many non-isomorphic  $p$-biased sets. It suffices to show that $f\mapsto S_f$ is injective. Let $S_f=S_g$. Equation \ref{Zdef} implies:
$$\bigsqcup_{(l,k)\in I}\bigsqcup_j X^l_{j+k,f}\times Y_j=\bigsqcup_{(l,k)\in I}\bigsqcup_j X^l_{j+k,g}\times Y_j$$
Consider $S_f\cap X^{l_0}\times Y_q=S_g\cap X^{l_0}\times Y_q$. Taking complements if $|\{k\in [q]: (l_0,k)\in I\}|=q-1$, this is equal to $X^{l_0}_{k,f}\times Y_q=X^{l_0}_{k,g}\times Y_q$ for some $k$, by the construction of $l_0$. Thus, $X^{l_0}_{k,f}=X^{l_0}_{k,g}$ and  $f=g$. Hence, $f\mapsto S_f$ is injective.\\
There are uncountably many sets $S_f$ and only countably many of those can be pairwise isomorphic, so there are uncountably many non-isomorphic  $p$-biased sets and functions.\blokje

\section{Open Problems}
The concept of a biased partition introducedin this paper raises the question what graphs contain them. Evidently a graph needs to $d$-regular for some $d$. A biased partition then consists of $d$ parts, which in turn consist of pairs of neighbors. Hence, we need $2d$ to divide the order of the graph. 
\begin{question}
What conditions on a graph imply the existence of a biased partition?
\end{question}

A more modest question towards a full characterization would be to identify a larger class of graphs allowing biased partitions. This paper shows their existence in $\mathbb{Z}^n$ and the paper by Gross and Grupel \cite{gross2017} implies the existence of biased a partition of the hypercube $\{0,1\}^n$ exactly if $n=2^k$ for some $k$.
\begin{question}
What discrete tori $\mathbb{Z}_{m}^n$ allow a biased partition?
\end{question}
Or more generally;
\begin{question}
What Cayley graphs allow a biased partition?
\end{question}

\section{Acknowledgement}
The author would like to thank B\'ela Bollob\'as for suggesting the problem and providing feedback on a draft of this paper.

\bibliography{bibliography}{}

\begin{thebibliography}{1}

\bibitem{benjamini1996}
Itai Benjamini and Harry Kesten.
\newblock Distinguishing sceneries by observing the scenery along a random walk
  path.
\newblock {\em Journal d'Analyse Math{\'e}matique}, 69(1):97--135, 1996.

\bibitem{finucane2014}
Hilary Finucane, Omer Tamuz, and Yariv Yaari.
\newblock Scenery reconstruction on finite abelian groups.
\newblock {\em Stochastic Processes and their Applications}, 124(8):2754--2770,
  2014.

\bibitem{gross2017}
Renan Gross and Uri Grupel.
\newblock Indistinguishable sceneries on the boolean hypercube.
\newblock {\em arXiv preprint arXiv:1701.07667}, 2017.

\bibitem{Howard1996}
C~Douglas Howard.
\newblock Detecting defects in periodic scenery by random walks on
  $\mathbb{Z}$.
\newblock {\em Random Structures \& Algorithms}, 8(1):59--74, 1996.

\bibitem{howard1996orthogonality}
C~Douglas Howard.
\newblock Orthogonality of measures induced by random walks with scenery.
\newblock {\em Combinatorics, Probability and Computing}, 5(3):247--256, 1996.

\bibitem{lindenstrauss1999}
Elon Lindenstrauss.
\newblock Indistinguishable sceneries.
\newblock {\em Random Structures and Algorithms}, 14(1):71--86, 1999.

\bibitem{matzinger2006}
Heinrich Matzinger and J{\"u}ri Lember.
\newblock Reconstruction of periodic sceneries seen along a random walk.
\newblock {\em Stochastic processes and their applications},
  116(11):1584--1599, 2006.

\end{thebibliography}
\bibliographystyle{plain}

\end{document}